\documentclass[11pt,a4paper]{article}
\usepackage{epsfig}
\usepackage{amsthm}

\usepackage[margin=1.1in]{geometry}
\usepackage{comment}
\usepackage{subfig}

\usepackage{amssymb,amsmath,latexsym,enumerate,graphicx,microtype,cite}
\usepackage{abstract}
\usepackage{url}
\usepackage{hyperref}
\allowdisplaybreaks
\usepackage{a4}
\usepackage{amsmath,amssymb,amsthm,mathrsfs}
\usepackage{mathtools}
\usepackage{comment}
\usepackage{abstract}
\usepackage{lineno}

\usepackage[colorinlistoftodos]{todonotes}

\renewcommand{\geq}{\geqslant}

\newcommand{\F}{{\cal{F}}} 
\newcommand{\A}{{\cal{A}}}

\newcommand{\bbox}{\vrule height7pt width4pt depth1pt}

\newtheorem{theorem}{Theorem}

\date{\today}

\author{Rom Pinchasi
\thanks{
Mathematics Department,
Technion---Israel Institute of Technology,
Haifa 32000, Israel.
{\tt room@technion.ac.il}. Visiting professor at EPFL, Lausanne, Switzerland. Supported by ISF grant (grant No.\ 1091/21)}
}

\title{A note on lenses in arrangements of pairwise intersecting circles in the plane} 

\begin{document}
\maketitle

\begin{abstract}
Let $\F$ be a family of $n$ pairwise intersecting circles in the
plane. We show that the number of lenses, that is convex digons,
in the arrangement induced by $\F$ is at most $2n-2$.
This bound is tight.
Furthermore, if no two circles in $\F$ touch, then 
the geometric graph $G$ on the set of centers of the circles in $\F$ whose edges correspond to the lenses generated by $\F$
does not contain pairs of avoiding edges. That is, $G$ does not
contain pairs of edges that are opposite edges in a convex quadrilateral.
Such graphs are known to have at most $2n-2$ edges.
\end{abstract}

\section{Introduction}
Given a family $\F$ of circles in the plane we consider 
the planar arrangement that is induced by the circles in
$\F$ and denote it by $\A(\F)$. The arrangement $\A(\F)$
consists of \emph{vertices} that are intersection points of circles in $\F$ and also of \emph{edges} that are arcs of circles in $\F$ delimited by two consecutive vertices.
Finally, $\A(\F)$ consists also of \emph{faces} that are the connected components of the plane after removing from it the union of all circles in $\F$.
A \emph{digon}
in $\A(\F)$ is a face in the arrangement $\A(\F)$ that has two edges. 
We distinguish between two types of digons in arrangements of circles.
A \emph{lens} in $\A(\F)$ is a face with two edges in $\A(\F)$ that is equal to the intersection of two discs bounded by circles in $\F$. Each of the
two circles in $\F$ corresponding to a lens is said to \emph{support} the lens. The lens is said to be 
\emph{created}
by these two circles supporting it.
Lenses are in fact just the convex digons in $\A(\F)$.
The arrangement $\A(\F)$ may contain also digons that are not convex. These are digons that are equal to the difference of two discs bounded by circles in $\F$. They are called \emph{lunes}.

Gr\"unbaum \cite{G72} conjectured that the number of digons in arrangements of $n$ pairwise intersecting pseudo-circles is at most $2n-2$.
This conjecture of Gr\"nbaum was verified in by Agarwal et al. \cite{ANPPSS04} for
arrangements of pseudo-circles surrounding a common point. In a recent work by Felsner, Roch, and Scheucher \cite{FRS23}, Gr\"unbaum's conjecture was verified for any arrangement of pairwise intersecting pseudo-circles under an additional assumption that the family of pseudo-circles contains three pseudo-circles every two of which create a digon in the arrangement.

In this paper we will be concerned with digons that are lenses in a family of pairwise intersecting circles. We will show that the number of lenses in  an arrangement of $n$ pairwise intersecting circles 
without any further assumption is at most $2n-2$.

\begin{theorem}\label{theorem:lenses}
Let $\F$ be a family of $n$ pairwise intersecting circles 
in the plane. Then $\A(\F)$ has at most $2n-2$ lenses.
This bound is tight for $n \geq 4$. 
\end{theorem}

\begin{figure}[ht]
	\centering
	\includegraphics[height=6cm]{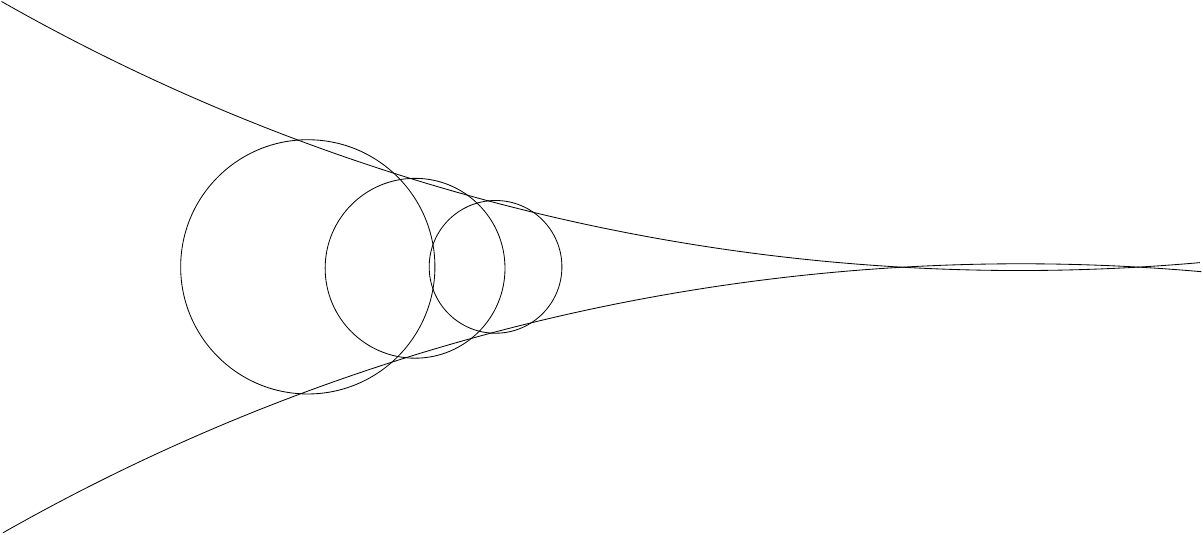}
	\caption{A family of $5$ pairwise intersecting circles with $8$ lenses.}
	\label{fig:tight}
\end{figure}

The simple construction in Figure \ref{fig:tight} shows 
that the bound in Theorem \ref{theorem:lenses} is best possible for $n \geq 4$. There are $5$ circles in this construction and $8$ lenses. One can generalize the construction for any number of circles by suitably   adding more circles to the three smaller circles
in the figure.

It is interesting to note that the number of lenses 
in arrangements of pairwise intersecting unit circles
is different. It is shown in \cite{P02} that 
there are at most $n$ lenses in any arrangement of
pairwise intersecting unit circles in the plane 
and this bound is best possible.

There has been a lot of research about lenses in arrangements of circles and pseudo-circles that are not necessarily pairwise intersecting. We will not survey here the vast literature about digons in arrangements of circles and pseudo-circles and on related situations where we allow curves to intersect more than twice and only refer the reader to \cite{GOT18} and the many references therein. The case where circles need not be pairwise intersecting is of completely different nature. We remark that in such a case the best constructions show that it is possible that $n$ circles will determine $\Omega(n^{4/3})$ many lenses. The best known upper bound is $O(n^{3/2}\log n)$ given in \cite{MT06}, that is following the footsteps of \cite{PR04}. The case of unit circles is of particular interest because of its relation to the celebrated unit distance problem posed by Paul Erd\H os \cite{E46}.    

Going back to families of pairwise intersecting circles, the number of lunes in these arrangements was studied in \cite{ALPS01}. 

\begin{theorem}[\cite{ALPS01}]\label{theorem:lunes}
Let $\F$ be a family of $n$ pairwise intersecting circles 
in the plane. Let $G$ be the geometric graph on the
set of centers of the circles in $\F$ whose edges 
correspond to pairs of discs bounded by circles in
$\F$ whose difference is a lune in $\A(\F)$.
Then $G$ is a bipartite planar embedding. Consequently,
$\A(\F)$ has at most $2n-4$ lunes.
\end{theorem}

Theorem \ref{theorem:lunes} is used in \cite{ALPS01}
to derive a linear upper bound (that is not tight) for the number of lenses
in arrangements of pairwise intersecting circles in the plane. Theorem \ref{theorem:lenses}, that we prove here, provides the tight upper bound for the number of lenses in a family of pairwise
intersecting circles in the plane.

Similar to the proof of Theorem \ref{theorem:lunes},
the proof of Theorem \ref{theorem:lenses} relies too
on studying the corresponding geometric graph $G=G(\F)$ whose vertices are the centers of circles in $\F$ and two centers are connected by an edge in $G(\F)$ iff the corresponding circles in $\F$ create a lens in $\A(\F)$.
We will show that unless we allow two circles in $\F$
to touch, then $G$ does not contain a pair of \emph{avoiding} edges. 
Two straight line segments (or edges in a geometric graph $G$) are called \emph{avoiding}
if they are opposite edges in a convex quadrilateral. Handling the case of $\F$ having pairs of touching circles
will require only a bit more effort. This is because if we allow touching circles in $\F$, the geometric graph $G(\F)$ may contain pairs of avoiding edges. Luckily, we will be able to compensate for this.

The important property of geometric graphs that do not contain pair of avoiding edges is given in \cite{KL98} and \cite{V98} (see also \cite{P08} for 
a different
and shorter proof of the same result, based on graph drawing). 

\begin{theorem}[\cite{KL98, V98}]\label{theorem:KLV}
Let $G$ be a geometric on $n$ vertices. If $G$ 
does not have a pair of avoiding edges, then it has at most $2n-2$ edges. This bound is tight for $n \geq 4$.
\end{theorem}

\bigskip

Our main goal is to prove Theorem \ref{theorem:lenses}.
In order to prove Theorem \ref{theorem:lenses} we start by making some assumptions without loss of generality that will help to simplify the presentation.

\begin{figure}[ht]
	\centering
	\includegraphics[width=6cm]{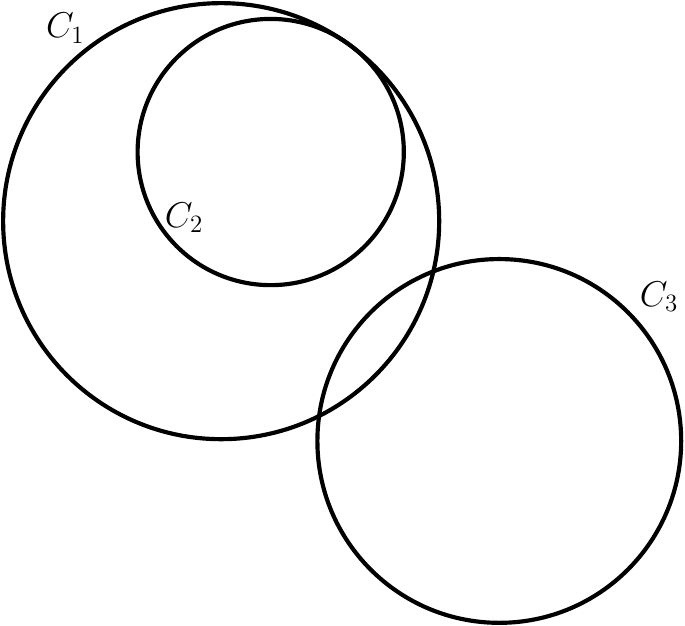}
	\caption{$C_{2}$ cannot be contained in the disc bounded by $C_{1}$.}
	\label{fig:not_contained}
\end{figure}

We first observe that we may assume that no two discs bounded by the circles in $\F$ may contain each other. 
Assume that $C_{1}$ and $C_{2}$ are two circles in 
$\F$ such that the disc bounded by $C_{1}$ fully contains
the disc bounded by $C_{2}$. Because every two circles 
in $\F$ intersect, it must be that $C_{2}$ touches $C_{1}$
internally. We claim that in such a case $C_{1}$ cannot
support any lens. The reason is that if $C_{1}$ create
a lens together with a circle $C_{3}$ in $\F$, then 
necessarily $C_{3}$ and $C_{2}$ cannot intersect (see Figure \ref{fig:not_contained}).

We can therefore remove $C_{1}$ from $\F$ and conclude 
Theorem \ref{theorem:lenses} by induction on the number 
of circles in $\F$.
We henceforth assume that no two discs
bounded by circles in $\F$ may contain each other. We conclude that if two circles in $\F$ touch, then 
they must touch externally and consequently no more
than two circles in $\F$ may be pairwise touching at a common point. 

\begin{figure}[ht]
	\centering
	\includegraphics[width=8cm]{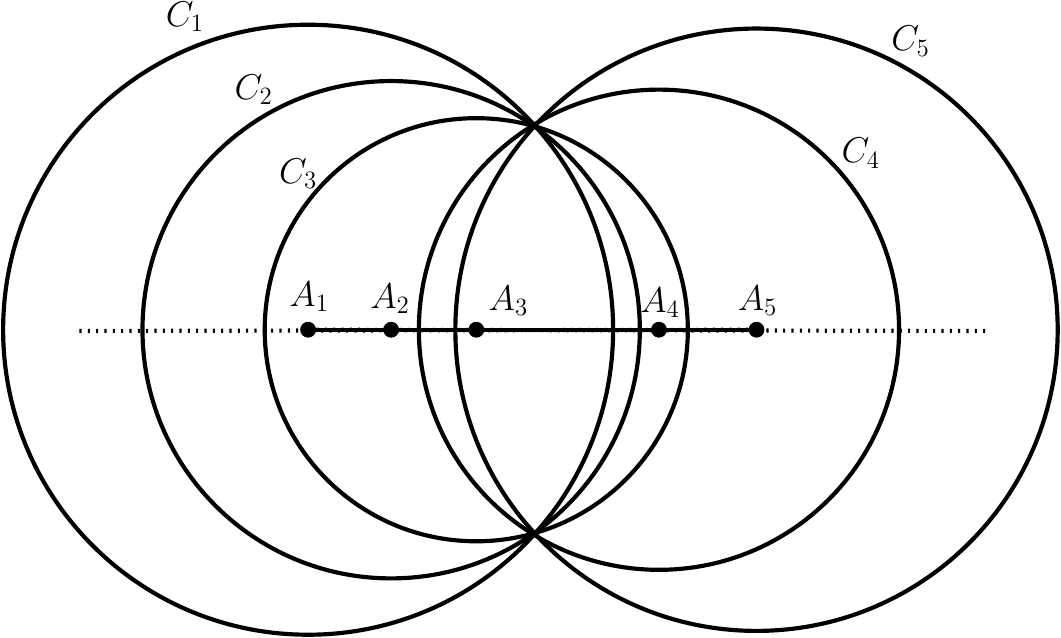}
	\caption{A pencil of circles. Only the two extreme circles, $C_{1}$ and $C_{5}$ in this figure, may create a lens.} \label{fig:pencil}
\end{figure}

For reasons that will become clear later, it will be more convenient for us to assume that
if $C_{1}$ and $C_{2}$ are two circles in $\F$ that support a lens, then the segment 
connecting the center of $C_{1}$ to the center of $C_{2}$ is not collinear with 
any center of circle in $\F$ that is not on that segment. We can indeed assume this without loss of generality. This is naturally the case if we assume that no three centers
of circles in $\F$ are collinear.
If we want to avoid this assumption, then 
we can just apply a generic inversion to the plane. In such a case three centers of circles in
$\F$ will remain collinear only if they 
form a \emph{pencil}, that is only if they
pass through two common points (or mutually touch at a point, which is impossible in our case). However, in such a case where we have 
several circles in $\F$ passing through two common points, then only the two extreme
circles may create together a lens (see Figure \ref{fig:pencil}).
In such a case the segment connecting the centers of these two extreme circles will not be collinear with any other center of a circle in $\F$ not on that segment. 
In order to prove Theorem \ref{theorem:lenses} and in order to allow also touching circles in $\F$, we define the geometric graph $G=G(\F)$
on the set of centers of the pairwise intersecting circles in $\F$ in the following way. The edges of $G$ will be either red, or blue.
We connect two centers of circles in $\F$ by a red edge if the corresponding
circles create a lens. We connect two centers of circles in $\F$ by
a blue edge if the corresponding circles touch.

The following theorem is
of independent interest and is also an intermediate step before proving Theorem 
\ref{theorem:lenses}.

\begin{theorem}\label{theorem:main}
The geometric graph $G(\F)$ does not contain pairs of avoiding edge unless
there are four circles $C_{1}, C_{2}, C_{3}$ and $C_{4}$ in $\F$
passing through a common point at which $C_{1}$ and $C_{3}$ touch each other and also $C_{2}$ and $C_{4}$
touch each other. Only in such a case it is possible to have avoiding edges in $G(\F)$. The avoiding edges in such a case can only be a pair of two opposite red edges in the quadrilateral whose vertices are the centers of such four circles $C_{1}, C_{2}, C_{3}$, and $C_{4}$.
\end{theorem}

We present the proof of Theorem \ref{theorem:main} in Section \ref{section:main}. Then we bring the proof of Theorem \ref{theorem:lenses} in Section \ref{section:lenses}.

\section{Proof of Theorem \ref{theorem:main}.}\label{section:main}
Assume that $e$ and $f$ are two avoiding
edges in $G(\F)$. Then $e$ and $f$
remain two avoiding edges in the graph $G$ that corresponds only to the four circles in $\F$ whose centers are the four endpoints of $e$ and $f$. This is because by
removing, or ignoring, all other circles in $\F$ we cannot destroy the digons, or pairs of touching circles corresponding to the edges $e$ and $f$. For the proof we will indeed assume that $\F$ is a family of only four circles. These are the four circles in the original family $\F$ that are centered at the endpoints of the avoiding segments $e$ and $f$. Consequetly, $G=G(\F)$ is a geometric graph with only four vertices. Denote by $A_{1}$ and $A_{2}$ the vertices of $e$
and let $A_{3}$ and $A_{4}$ be the vertices of $f$.
Because $e$ and $f$ are avoiding, we assume without loss of generality that 
$A_{1}A_{2}A_{3}A_{4}$ is a convex quadrilateral.
For $i=1,2,3,4$ denote by $C_{i}$ the circle in 
$\F$ centered at $A_{i}$. Denote by $D_{i}$ the
closed circular disc bounded by $C_{i}$.

In order to simplify the presentation of the proof we would like
to assume that there is a point $M$ common to three of the circles $C_{1}, C_{2}, C_{3}$, and $C_{4}$. We can assume this without loss of generality. This is because if this is not the case we  inflate the circle $C_{1}$ keeping its center $A_{1}$ fixed until the first time $C_{1}$ passes through an intersection point of two other circles from $\F$, namely an intersection point of two of the circles $C_{2}, C_{3}$, and $C_{4}$.

More precisely, assume that no three of $C_{1}, C_{2}, C_{3}$, and $C_{4}$
pass through a common intersection point. 
Because $C_{1}$ and $C_{2}$ create a lens or externally touch, it must be that the intersection points on $C_{2}$ with the circles $C_{3}$ and $C_{4}$ lie outside of $D_{1}$.
We start inflating $C_{1}$  until the first time it meets an intersection points of two of $C_{2}, C_{3}$, and $C_{4}$.
When this happens $C_{1}$ still intersects with each of 
$C_{2}, C_{3}$, and $C_{4}$
because each of $C_{2}, C_{3}$, and $C_{4}$ contains at least one intersection point that is not surrounded by $C_{1}$, 
while each of $C_{2}, C_{3}$, and $C_{4}$ contains also points in $D_{1}$ because before we inflated $C_{1}$ it intersected each of $C_{1}, C_{2}$, and $C_{3}$, while $D_{1}$ only increases through the inflation of $C_{1}$.
The only thing that is left two show is that $e$ and $f$ are still edges in $G(\F)$, that is, we need to show that even after the inflation of $C_{1}$ it is still true that $C_{3}$ and $C_{4}$ create a lens or touch and the same is true for $C_{1}$ and $C_{2}$.
The reason this is true is that the intersection points of $C_{1}$
with any of $C_{2}, C_{3}$, and $C_{4}$ move continuously with the inflation of $C_{1}$. Therefore, if $C_{3}$ and $C_{4}$ create a lens, then the edges $C_{3} \cap D_{4}$ and $C_{4} \cap D_{3}$, remain edges in $\A(\F)$ also after inflating $C_{1}$. It is also clear that if $C_{3}$ and $C_{4}$ touch, then they remain touching regardless of the inflation of $C_{1}$. If $C_{1}$ and $C_{2}$ create a lens, then 
the edges $C_{1} \cap D_{2}$ and $C_{2} \cap D_{1}$ remain edges in $\A(\F)$ also after inflating $C_{1}$
and therefore $C_{1}$ and $C_{2}$ create a lens also after we inflate $C_{1}$.

The only case where we need to be careful is if $C_{1}$ and $C_{2}$ touch each other
before we inflate $C_{1}$.
In this case as we start inflating $C_{1}$ it creates a lens with $C_{2}$ and by the
argument above it will create a lens with $C_{2}$ also when
the inflation of $C_{1}$ is stopped.

We conclude that by inflating 
$C_{1}$,
keeping its center fixed, we may assume that three of the circles in $\F$ pass through a common point while $e$ and $f$ remain two avoiding edges in $G(\F)$. Therefore, we assume without loss of generality that $C_{1}, C_{2}$, and $C_{3}$ pass through a common intersection point that we denote by $M$. 

\begin{figure}[ht]
	\centering
	\includegraphics[width=8cm]{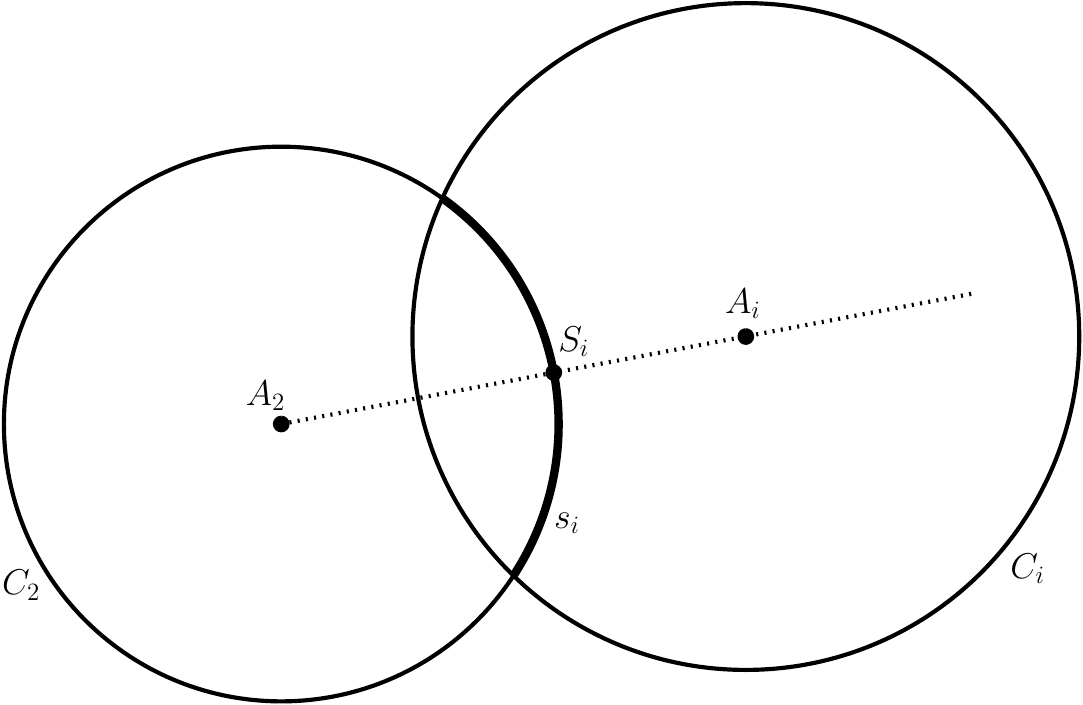}
	\caption{The arc $s_{i}$ with its center $S_{i}$.}	\label{fig:s_i}
\end{figure}

\begin{figure}[ht]
	\centering
	\includegraphics[width=5.1cm]{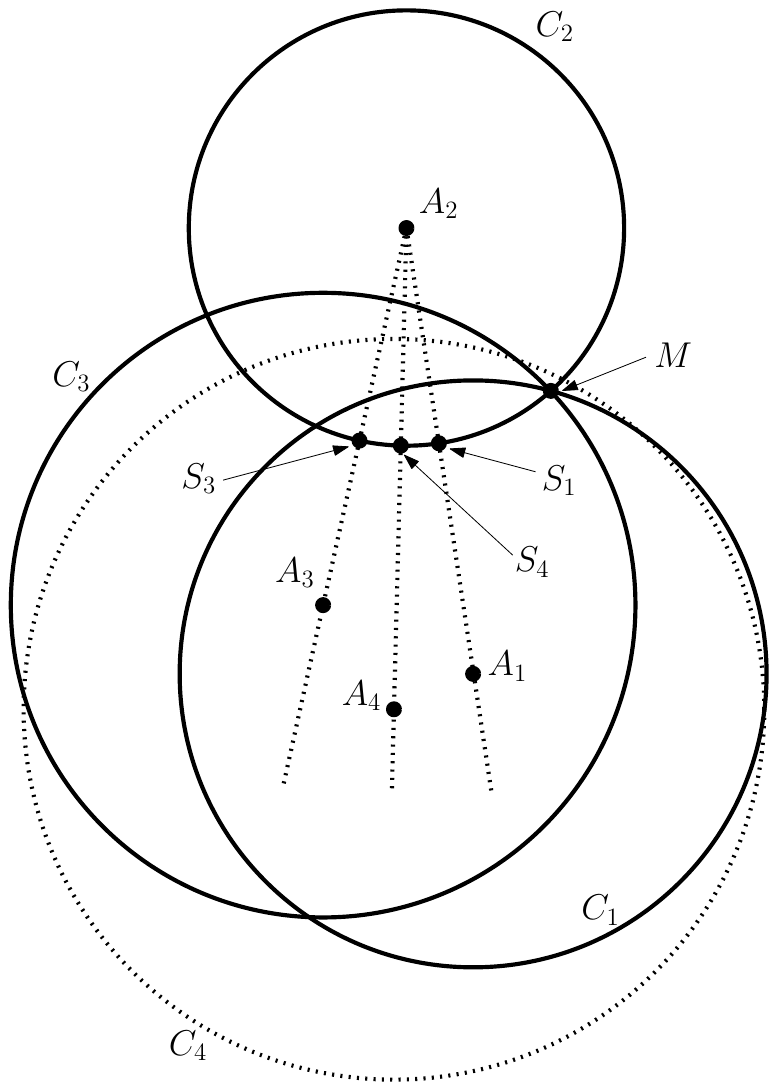}
	\caption{The impossible case where $s_{3}$ contains $s_{1}$.}
	\label{fig:s_1_s_3}
\end{figure}

Denote by 
$s_{1},s_{3}$, and $s_{4}$ the arcs 
$C_{2} \cap D_{1}, C_{2} \cap D_{3}$, and $C_{2} \cap D_{4}$, respectively. We note that each $s_{i}$ may be degenerate and equal to a single point in case 
$C_{i}$ and $C_{2}$ touch.
Notice that $M$ is an endpoint of both 
$s_{1}$ and $s_{3}$. For $i=1,3,4$ let $S_{i}$ denote the center of the arc $s_{i}$.
We notice that $S_{i}$ is the point of intersection of the ray $\overrightarrow{A_{2}A_{i}}$ with $C_{2}$. 
For this reason and because $A_{1}A_{2}A_{3}A_{4}$ is a convex quadrilateral, the point $S_{4}$ must lie in the shorter
arc of $C_{2}$ delimited by $S_{1}$ and $S_{3}$ (see Figures \ref{fig:s_i} and \ref{fig:s_1_s_3}).

We claim that the intersection of $s_{1}$ and $s_{3}$ is equal to the point $M$.
Indeed, otherwise $s_{1}$ and $s_{3}$ are both nondegenerate arcs and the relative interiors
of $s_{1}$ and $s_{3}$ overlap.
Then arc $s_{1}$ is an edge of the lens $D_{1} \cap D_{2}$ and cannot contain any intersection points in its relative interior.
The arc $s_{1}$ can also not be equal to $s_{3}$, or else $A_{1}, A_{2}$, and $A_{3}$ are collinear, contrary to the assumption that they are three vertices of a convex quadrilateral. We conlude
from here that it must be that $s_{3}$ strictly contains $s_{1}$ (see Figure \ref{fig:s_1_s_3}).

Because $S_{4}$ lies in the shorter arc of $C_{2}$ delimited by $S_{1}$ and $S_{3}$, it follows that $S_{4}$ lies in the relative interior of $s_{3}$. Consequently, $S_{4}$ lies in the interior of $D_{3}$. At the same time $S_{4}$ lies also in $D_{4}$. This shows that $S_{4}$, that is a point on $C_{2}$, belongs to the lens or is equal to the touching point $D_{3} \cap D_{4}$.
This is possible only if $S_{4}$ is 
an intersection point of $C_{3}$ and $C_{4}$. However, the latter case is impossible because $S_{4}$ lies in the interior of $D_{3}$.

\begin{figure}[ht]
	\centering
	\includegraphics[width=6cm]{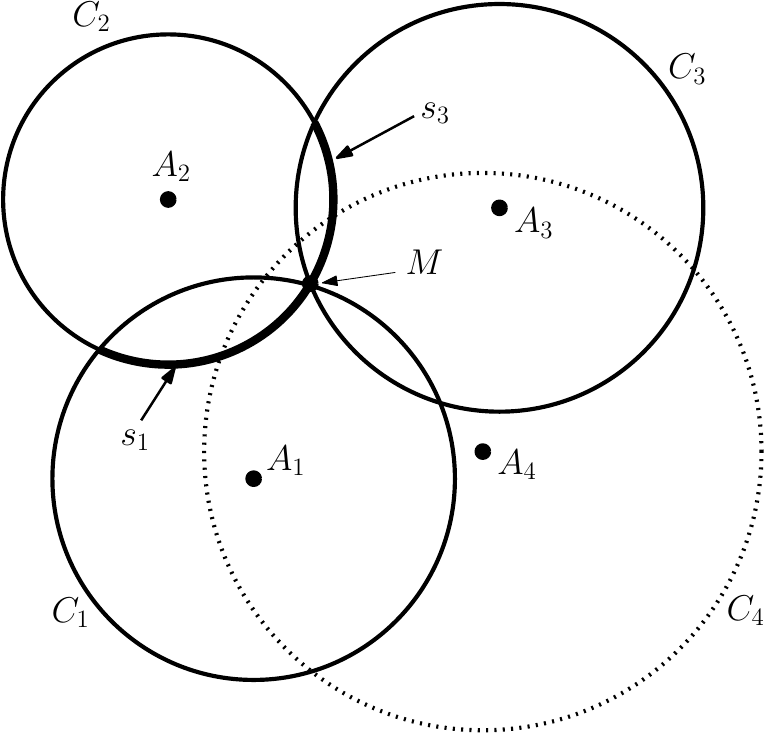}
	\caption{$M$ cannot belong to the relative interior of $s_{4}$.}
 \label{fig:M_not_intr}
\end{figure}

Having shown that the intersection of $s_{1}$ and $s_{3}$ is equal to the point $M$,
we observe that
$M$ cannot belong to the interior of
$D_{4}$, or else $D_{3} \cap D_{4}$ cannot be a lens nor a touching point. Consequently, 
$M$ cannot belong to the relative interior of the arc $s_{4}$.
Combining this with the fact that $S_{4}$
belongs to the shorter arc of $C_{2}$ delimited by $S_{1}$ and $S_{3}$, we conclude that $s_{4}$ is contained in $s_{1}$, or it is contained in $s_{3}$ (see Figure \ref{fig:M_not_intr}).

We claim that $s_{4}$ must be equal to the point $M$ and hence $C_{2}$ and $C_{4}$ touch each other at $M$.
In order to prove this we consider the two possible cases
$s_{4} \subset s_{1}$, and 
$s_{4} \subset s_{3}$.

\medskip

\begin{figure}[ht]
	\centering
	\includegraphics[width=8cm]{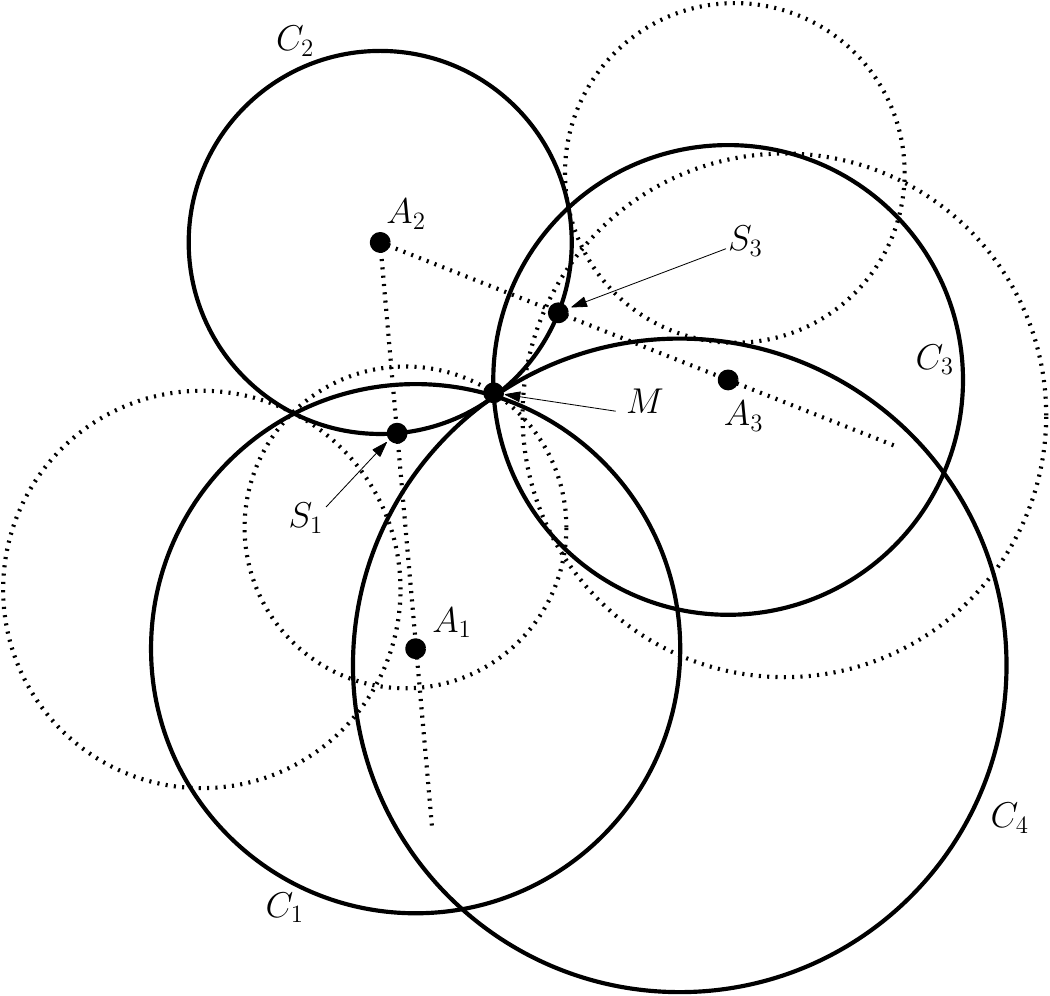}
	\caption{$C_{4}$ must touch $C_{2}$ at $M$. The dotted circles show impossible positions of $C_{4}$.}	\label{fig:C_4_touch_M}
\end{figure}

\noindent {\bf Case 1.}
$s_{4} \subset s_{1}$. In this case $s_{1}$ must be a nondegenerate arc, or else 
both $C_{1}$ and $C_{4}$ 
touch $C_{2}$ at the same point $M$ which is impossible.
Therefore, $s_{1}$ is an edge of the lens $D_{1} \cap D_{2}$ and it cannot contain any intersection points in its relative interior. Hence
either $s_{4}$ is equal to an endpoint of $s_{1}$, or $s_{4}$ is equal to $s_{1}$. The latter case is impossible as it would imply that $A_{1}, A_{2}$, and $A_{4}$ are collinear, contrary to the assumption that they are vertices of a convex quadrilateral. In the former case $s_{4}$ must be equal to the point $M$ because $S_{4}$,
that is equal to $s_{4}$ in this case, lies in the shorter arc of $C_{2}$ delimited by $S_{1}$ and $S_{3}$ (see Figure \ref{fig:C_4_touch_M}).

\medskip

\noindent {\bf Case 2.} 
$s_{4} \subset s_{3}$. In this case $s_{3}$ must be a nondegenerate arc, or else 
both $C_{1}$ and $C_{3}$ 
touch $C_{2}$ at the same point $M$ which is impossible. 
We notice that $s_{4}$ is contained in both $D_{4}$ and $D_{3}$ and therefore $s_{4}$ is contained in the lens $D_{4} \cap D_{3}$. 
This is impossible unless $s_{4}$ is a single point
that is equal to an intersection point of $C_{3}$ and $C_{4}$ and in particular belongs to $C_{3}$. Then $s_{4}$ must be equal to one of the two endpoints of $s_{3}$. It must be equal to the point $M$ because $S_{4}$,
that is equal to $s_{4}$ in this case, lies in the shorter arc of $C_{2}$ delimited by $S_{1}$ and $S_{3}$ (see Figure \ref{fig:C_4_touch_M}).

\begin{figure}[ht]
	\centering
	\includegraphics[width=6cm]{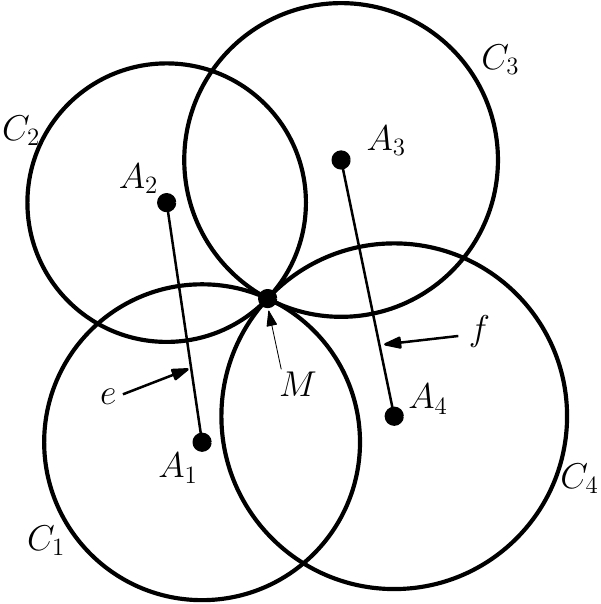}
	\caption{The only configuration in which
 $e$ and $f$ may be avoiding. Both pairs of circles $C_{2}$ and $C_{4}$ as well as $C_{1}$ and $C_{3}$ mutually touch at $M$.}\label{fig:special_conf}
\end{figure}

Having shown that that $C_{1}, C_{2}$, and $C_{4}$ pass through $M$ we can argue  symmetrically that $C_{3}$ touches $C_{1}$
at $M$, as illustrated in Figure \ref{fig:special_conf}.

We observe that in such a case the circle $C_{1}$ could not be inflated by a factor greater than $1$, or else it would be disjoint from $C_{3}$ in the original configuration of the circles.
To summarize, we have shown that if $e$ and $f$ are two avoiding edges in $G(\F)$, then 
necessarily $C_{1}$ touches $C_{3}$ at a point $M$ and $C_{3}$ touches $C_{4}$ at the same point $M$. This concludes the proof of Theorem \ref{theorem:main}.
\bbox

\section{Proof of Theorem \ref{theorem:lenses}.}\label{section:lenses}

We need to show that $\F$ determines at most $2n-2$ lenses.
We consider the geometric graph $G=G(\F)$ defined 
just before the statement of Theorem \ref{theorem:main}. Proving Theorem
\ref{theorem:lenses} is equivalent to showing that the number of red edges in $G(\F)$
is at most $2n-2$.

\begin{figure}[ht]
	\centering
	\includegraphics[width=10cm]{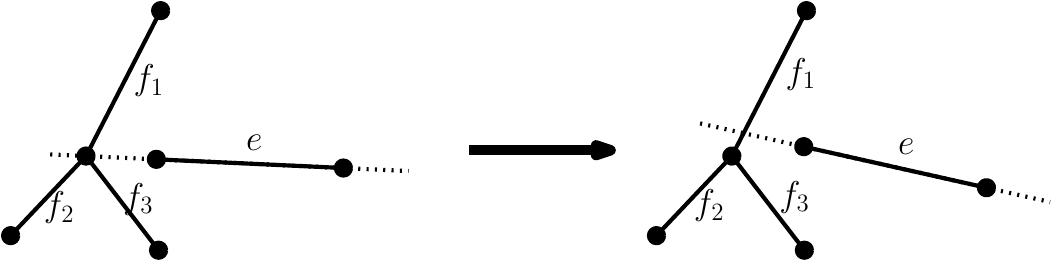}
	\caption{The edge $e$ is collinear with a vertex not on $e$.}\label{fig:perturb}
\end{figure}

If no two circles in $\F$ touch each
other, then it follows from Theorem \ref{theorem:main} that the graph $G(\F)$ consists of red edges only. Moreover, 
no two edges in $G(\F)$ are avoiding.
If the set of vertices of $G(\F)$, that is the set of centers of the circles in $\F$,
is in general position in the sense that no three
of them are collinear, then we can apply Theorem \ref{theorem:KLV}
and conclude that $G(\F)$ has at most $2n-2$ edges. Consequently, $\F$ determines at most 
$2n-2$ lenses. If the set of vertices of $G(\F)$
is not in general position,
then strictly speaking we cannot apply Theorem \ref{theorem:KLV} because Theorem \ref{theorem:KLV}, as most other theorems about geometric graphs, is stated and proved for geometric graphs whose set of vertices is in general position. This is a standard assumption in most results concerning geometric graphs.
In order to be able to apply Theorem \ref{theorem:KLV}, we perturb a bit the vertices of $G(\F)$
to make them lie in general position. We notice that by perturbing the vertices of $G(\F)$
we may not create pairs of avoiding edges, 
unless there is an edge in $G(\F)$ that is collinear with another vertex in $G(\F)$ not on this edge (see Figure \ref{fig:perturb}, where in the perturbed picture on the right $e$ and $f_{2}$ are avoiding as well as $e$ and $f_{3}$). 

As we assumed without loss of generality, this cannot not happen. No edge $e$ in $G(\F)$ is collinear with a vertex not on on $e$.
We may therefore apply Theorem \ref{theorem:KLV} to
the perturbed geometric graph and conclude that $G(\F)$ has at most $2n-2$ edges as before. This proves Theorem \ref{theorem:lenses} in the case where no pair of circles in $\F$ touch each other.

In the more general case where we allow circles in $\F$ to touch
each other, the graph $G(\F)$ may contain pairs of avoiding red edges. However, in this case we can use the blue edges in $G(\F)$ as follows.
By Theorem \ref{theorem:main}, whenever the graph $G(\F)$ contains pairs of avoiding edges
it must be because of the very special structure as described in the statement
of Theorem \ref{theorem:main} and shown in Figure \ref{fig:special_conf}. If $e$ and $f$ are two
avoiding edges in $G(\F)$, then they are necessarily red edges that correspond to two lenses in 
$\A(\F)$. Let $A_{1}A_{2}A_{3}A_{4}$ be the convex quadrilateral
such that $e=A_{1}A_{2}$ and $f=A_{3}A_{4}$.
We assume without loss of generality that $A_{1}, A_{2}, A_{3}$, and $A_{4}$ is the clockwise cyclic order of
these points as vertices of the convex quadrilateral $A_{1}A_{2}A_{3}A_{4}$ (see Figure \ref{fig:special_conf}).
For $i=1,2,3,4$ denote by $C_{i}$ the circle
in $\F$ centered at $A_{i}$. By Theorem \ref{theorem:main}, $A_{1}$ and $A_{3}$ are connected by a blue edge in $G(\F)$ and so are $A_{2}$ and $A_{4}$. This corresponds to that $C_{1}$ and $C_{3}$ must touch, at a point that we denote by $M$, and so are $C_{2}$ and $C_{4}$. 

\begin{figure}[ht]
	\centering
	\includegraphics[width=9cm]{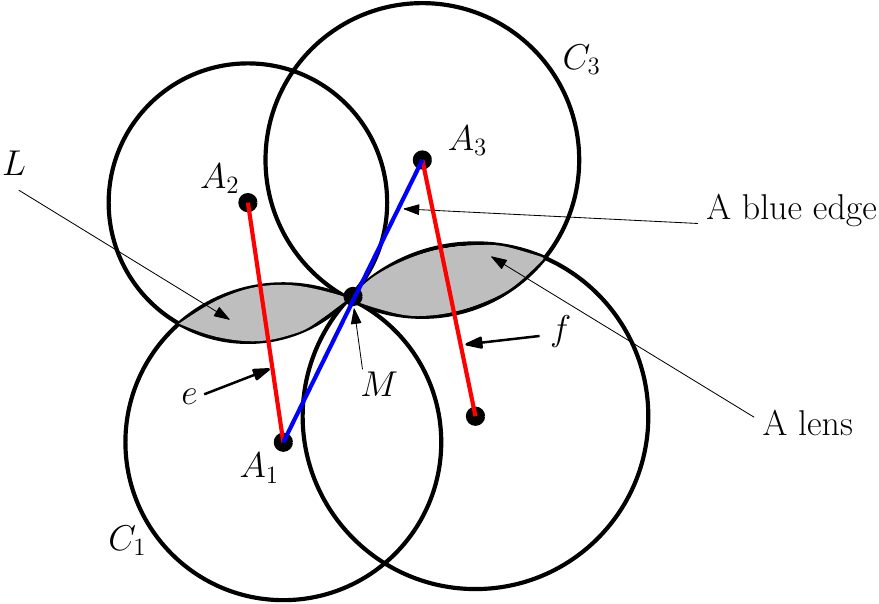}
	\caption{The blue edge $A_{1}A_{3}$ is uniquely charged.}\label{fig:special_charge}
\end{figure}

We remove either the red edge $e$, or the red edge $f$ 
from $G(\F)$ and charge the avoiding pair $e$ and $f$ to the 
blue edge $A_{1}A_{3}$. We claim that the blue edge
$A_{1}A_{3}$ cannot be charged for another red edge in this way.
This is because once we fix $A_{1}$ and $A_{3}$ and therefore
also $C_{1}$ and $C_{3}$, then we determine also the touching point
$M$ of $C_{1}$ and $C_{3}$. We claim that 
the edges $e$ and $f$ are determined as well. Indeed, $e$ is the edge that correspond to the unique lens $L$ supported by $C_{1}$ 
such that $M$ is a vertex of $L$ in $\A(\F)$ and $A_{1},A_{2}, M$ is the clockwise order of these three points, where $A_{2}$ is the center of the other circle in $\F$ supporting $L$ (see Figure \ref{fig:special_charge}). Any other such lens $L'$ would overlap with $L$, which is impossible. 

By a symmetric argument we show that the edge $f$ is determined by the blue edge $A_{1}A_{3}$. We conclude from this that the blue edge $A_{1}A_{3}$ in $G(\F)$ can be charged to at most one pair of avoiding edges in $G(\F)$. 

After repeating this procedure for every remaining pair of 
avoiding edges we are left with a subgraph $G'$ of $G(\F)$ 
in which no two edges are avoiding. The number of edges, red and blue, in $G'$ 
is greater than or equal to the number of red edges in $G$.
We apply Theorem \ref{theorem:KLV} to $G'$,
after possibly perturbing the vertices of $G'$, as we did already before, and conclude that $G'$ has at most $2n-2$ 
edges.
Therefore, there are at most $2n-2$ red edges in $G(\F)$. Consequently, there are at most $2n-2$ lenses in $\A(\F)$, as desired.
\bbox

\bigskip

\noindent {\bf \Large Acknowledgements.}
We thank Eyal Ackerman for 
helpful comments about the case where the centers of the circles in $\F$ are not in general position.

\bibliographystyle{abbrv}

\end{document}